\documentclass{article}
\usepackage{latexsym}
\begin{document}
%%%%%%%%%%%%%%%%%%%%%%%%%%%%%%%%%%%%%%%%%%%%%%%%%%%%%%%%%%%%%%%%%%%%%%%%%%%%%%%
%%%%%%%%%%%%%%%%%%%%%%%%%%%%%%%%%%%%%%%%%%%%%%%%%%%%%%%%%%%%%%%%%%%%%%%%%%%%%%%
%%%%%%%%%%%%%\hfill S. A. Skarst Choroszavin. \hfill \\

%\input unitary-2013fev20-0.tex

\begin{center}
{\bf  An Example of $J$-unitary Operator.
 Solving a Problem Stated by M.G. Krein.}
\bigskip

 {{\sc Sergej A. Choroszavin } }

 {\it mailto: choroszavin@narod.ru }
\bigskip

\bigskip

\end{center}

\begin{abstract}
{\bf Theorem 1.}
 Given a number $c \geq 1$, there exists a $J$-unitary operator $\widehat{ V}$, such that:  
\begin{eqnarray*} 
 &(a)& \quad  r(\widehat{ V} )= r(\widehat{ V}^{-1} )=  c ; 
\\
 &(b)&
S(\widehat{V},c)=S(\widehat{V}^{-1},c)
=S(\widehat{V}^{*-1},c) = S(\widehat{V}^*,c)=\{0\} \,.\,
\\
&(c)&\mbox { there exist maximal strictly positive and strictly negative } 
 \\                             
&&\widehat{ V}^{\pm}-\mbox{invariant subspaces }
\\
&& {\cal L_+, L_- } \mbox{ such that they are mutually $J$-orthogonal and } 
\\
&& {\cal L_++ L_- } \mbox { is dense in the space. }
\\
&(d_1)&\mbox { if $L_1$ is non-zero $\widehat{ V}$-invariant subspace, }
\widehat{ V}L_1 \subset L_1 \,, 
\\
&&
\mbox{ then } r(\widehat{ V}|L_1)=r(\widehat{ V}) \,.    
 \\                             
&(d_2)&\mbox { if $L_2$ is non-zero $\widehat{ V}^{-1}$-invariant subspace, }
\widehat{ V}^{-1}L_2 \subset L_2 \,, 
\\
&&
\mbox{ then } r(\widehat{ V}^{-1}|L_2)=r(\widehat{ V}^{-1}) \,.    
 \\                             
\end{eqnarray*}
\end{abstract}

In this paper we try to essencially refine some results of  
0404071.
So, the paper inherit the system of notations and definitions of 0404071.
For the references, historical comments, etc.,  see 0404071 as well.

\par\addvspace{\bigskipamount}\par\noindent
To specify the Problem, 
recall some definitions and facts

\par\addvspace{\bigskipamount}\par\noindent
{\bf Definition D1-1}
Given an operator $T$  and a number $c$, we write 
$$
 S(T,c)\,:=\,\{x\in H|\, \exists M\geq 0 \, \forall N\geq 0 
 \quad\|T^Nx\|\leq Mc^N \},
$$ 
$$
 r(T):= \hbox{spectral radius of }T \,.
$$

\par\addvspace{\bigskipamount}\par\noindent
{\bf Remark R1-1}
 Given a linear bounded operator 
$T$
 and a 
$T$-invariant subspace
$L$
 such that 
$r(T|L)\leq c$, 
 then 
$$
 L\subset S(T, c+\epsilon)
$$
for any $\epsilon > 0$.

\par\addvspace{\bigskipamount}\par\noindent
{\bf Remark R1-2}
$$
 S(T_1\oplus T_2,c)=S(T_1,c)\oplus S(T_2,c);
$$
 
\par\addvspace{2\medskipamount}\par\noindent 
The central  theorem of this paper is 
\par\addvspace{2\medskipamount}\par\noindent 
{\bf Theorem 1.}
 Given a number $c \geq 1$, there exists a $J$-unitary operator $\widehat{ V}$, such that:  
\begin{eqnarray*} 
 &(a)& \quad  r(\widehat{ V} )= r(\widehat{ V}^{-1} )=  c ; 
\\
 &(b)&
S(\widehat{V},c)=S(\widehat{V}^{-1},c)
=S(\widehat{V}^{*-1},c) = S(\widehat{V}^*,c)=\{0\} \,.\,
\\
 &(c)&\mbox { there exist maximal strictly positive and strictly negative } 
 \\                             
&&\widehat{ V}^{\pm}-\mbox{invariant subspaces }
\\
&& {\cal L_+, L_- } \mbox{ such that they are mutually $J$-orthogonal and }
\\
&& {\cal L_++ L_- } \mbox { is dense in the space. }
\\
&(d_1)&\mbox { if $L_1$ is non-zero $\widehat{ V}$-invariant subspace, }
\widehat{ V}L_1 \subset L_1 \,, 
\\
&&
\mbox{ then } r(\widehat{ V}|L_1)=r(\widehat{ V}) \,.    
 \\                             
&(d_2)&\mbox { if $L_2$ is non-zero $\widehat{ V}^{-1}$-invariant subspace, }
\widehat{ V}^{-1}L_2 \subset L_2 \,, 
\\
&&
\mbox{ then } r(\widehat{ V}^{-1}|L_2)=r(\widehat{ V}^{-1}) \,.    
 \\                             
\end{eqnarray*}
\par\addvspace{2\medskipamount}\par\noindent 
{\bf Proof .}
 By the Lemma 6 below there exists an operator $V$ in a separable 
 (real or complex, as desired) Hilbert space, such that
 $V$ is bounded, $V^{-1}$ exists, is bounded, and  
 $$ 
r( V )= r({ V}^{-1} )=  c
$$
$$
S(V,c)=S(V^{-1},c)=S(V^{*-1},c) = S(V^*,c)=\{0\} \,.\,
$$
Now put 
$$ \widehat{ V} = V \oplus {{V}^{*}}^{-1} $$
 and the proof of (a), (b) is evident.
To prove (c) apply Lemma 7. 
\\
Now proof of $(d_1)$.
\\
Suppose 
$r(\widehat{ V}|L_1)<r(\widehat{ V})$.
Put $\epsilon=r(\widehat{ V})-r(\widehat{ V}|L_1)$. Then $\epsilon>0$ and   
$$
 L_1\subset S(\widehat{ V}, r(\widehat{ V}|L_1)+\epsilon)
= S(\widehat{ V},r(\widehat{ V})) =\{0\} \,.
$$
Contradiction with $L_1\not=\{0\}$.
Proof of $(d_2)$ is quite similar.
\\
$\Box$

\par\addvspace{2\medskipamount}\par\noindent 
Now the details.
\par\addvspace{2\medskipamount}\par\noindent 

Let 
$\{u_n\}_{n\in {\bf Z}}$ 
 be a bilateral number sequence; we will suppose that 
$u_n \not= 0$
 for all 
$n\in {\bf Z}$.

 In this case let 
$U$ 
 denote the {\bf shift } 
\footnote{ 
 the full name is: the {\bf bilateral weighted shift of 
$\{b_n\}_n$,
 to the right}.
} %%%% end of footnote 
 that is generated by the formula 
$$
 U\,:\,b_n\,\mapsto\,\frac{u_{n+1}}{u_n}\,b_{n+1} \qquad . \eqno(*)
$$

 The general facts we need are these:

{\bf Observation O2-1}

 One constructs the 
$U$ as follows:
 
 One starts extending the instruction
$(*)$
 on the linear span of the 
$\{b_n\}_{n\in {\bf Z}}$ 
 so that the resulted operator becomes linear.
 That extension is unique and defines a linear densely defined 
 operator, which is here denoted by
$U_{min}$,
 and which is closable.
 The closure of
$U_{min}$
 is just the 
$U$.

 Now then, this 
$U$
 is closed and at least densely defined and injective;
 it has dense range and the action of
$U^N$, 
$U^{*-N}$, 
$U^{*N}U^N$, 
$U^{-N}U^{*-N}$ 
 (for any integer $N$) is generated by 
$$
 U^N:b_n\mapsto \frac{u_{n+N}}{u_n}b_{n+N}\,;\qquad 
 U^{*-N}:b_n\mapsto \frac{u_n^*}{u_{n+N}^*}b_{n+N}\,;
$$
$$
 U^{*N}U^N:b_n\mapsto {|\frac{u_{n+N}}{u_n}|}^2 b_n\,;\qquad 
 U^{-N}U^{*-N}:b_n\mapsto {|\frac{u_n}{u_{n+N}}|}^2 b_n\,.
$$

 In particular,
$U^N$
 is bounded just when the number sequence 
$\{|u_{n+N}/u_n|\}_n$
 is bounded.
 Moreover,
$$  \| U^N \| = sup\{ \; |u_{n+N}/u_n| \;  | n \in {\bf Z} \}    $$ 
\\ 
$\Box$

 The special factors we need are : 

{\bf Observation O2-2}(revised)
Let $a>0$ be number.  Then
$$
\begin{array}{llcccccccccc}

S(U,a)\not= \{0\}          & \Leftrightarrow && \exists M^{'} 
 \,    |u_N|         &    \leq      &   M^{'}a^{N}    & 
 \,  (N=0,1,2,...)   &  \\

S(U^{*-1},a)\not= \{0\}         & \Leftrightarrow &&\exists M^{'}
 \,    |u_N|^{-1}    &     \leq     &   M^{'}a^{N}    &
 \,  (N=0,1,2,...)   &  \\

 S(U^{-1},a)\not= \{0\}               & \Leftrightarrow &&\exists M^{'}
 \,    |u_{-N}|      &     \leq     &    M^{'}a^{N}   &
 \,  (N=0,1,2,...)   &   \\

 S(U^{*},a) \not= \{0\}         & \Leftrightarrow && \exists M^{'}
 \,    |u_{-N}|^{-1} &     \leq     &    M^{'}a^{N}   &
 \,    (N=0,1,2,...) &     \\
\end{array}
$$
{\bf Proof.}
 The family 
$\{b_n\}_n$ 
 is an orthonormal basis. In addition 
$U^Nb_n\perp U^Nb_m$ for 
$n\neq m$. Thus 
$$
 {\|U^Nf\|}^2=\sum_n|(b_n,f)|^2\|U^Nb_n\|^2=
 \sum_n|(b_n,f)|^2{|\frac{u_{n+N}}{u_n}|}^2\,
$$
 for every 
$f\in D_{U^N}$.

 In particular, 
$$\|U^Nf\|\geq|(b_n,f)||u_{n+N}/u_n|$$ 
 for all integers 
$n$. 

 It follows that: 

 Given 
$f\in H_0\setminus \{0\}$
 and given some real $M,a$ such that 
$$
 \|U^Nf\|\leq Ma^N 
 \mbox{ for } 
 N=0,1,2,...\,,
$$
 then there exists a real 
$M'$
 such that 
$$
 |u_N|\leq M'a^N 
 \mbox{ for } 
 N=0,1,2,...
$$

On the other hand, if
$$
 \exists M^{'} 
 \,    |u_N|            \leq        M^{'}a^{N}    
 \,  (N=0,1,2,...)  
$$
 then  $b_0 \in S(U,a)$ and $S(U,a)\not=\{0\}$. 
 
For
$U^{*-1},U^{-1},U^*$, the proof is similar. 
\par\addvspace{2\medskipamount}\par\noindent
\par\noindent

\newpage

\par\addvspace{2\medskipamount}\par\noindent
\par\noindent 
% %
{\bf Lemma 1.}
\par\noindent 
Let
\begin{eqnarray*}
|v_n|  
&\not=&
0
\end{eqnarray*}
and $V$ be a bilateral weighted shift defined by 
\begin{eqnarray*}
Vb_n  
&:=&
\frac{v_{n+1}}{v_n} b_{n+1}
\end{eqnarray*}
Suppose that
$$
      \exists n_0 \geq 0  
$$
so that
$$
  \forall n  \qquad 
|n| \geq n_0  \Rightarrow    \Bigl| \frac{v_{n+1}}{v_n} \Bigr| \leq  1   
\qquad\qquad\qquad\qquad(*)
$$
 Then
$$||V||<\infty$$
and 
\begin{eqnarray*}
 r(V)
&\leq&
1 
\end{eqnarray*}
As a consequence 
\begin{eqnarray*}
 r(V^*)=r(V)
&\leq&
1 
\end{eqnarray*}
\medskip\noindent
\par\addvspace{2\medskipamount}\par\noindent
{\bf Proof.}
\\
By supposition,
$$
|n| \geq n_0  \Rightarrow    \Bigl| \frac{v_{n+1}}{v_n} \Bigr| \leq  1   
\qquad\qquad\qquad\qquad(*)
$$
Put 
$$
  c_0 = max\{ 1, max \Bigl\{ \Bigl| \frac{v_{n+1}}{v_n} \Bigr| \quad 
             \Big| \quad n=-n_0, -n_0+1, \cdots , n_0  \Bigr\}\}
$$
Hence for every $n \in {\bf Z}$, $N > 0$ 
$$
  \Bigl| \frac{v_{n+N}}{v_n} \Bigr| = \Bigl| \frac{v_{n+1}}{v_n} \Bigr|\cdots\Bigl| \frac{v_{n+N}}{v_{n+N-1}} \Bigr|  \leq {c_0}^{2n_0+1} 
$$
Apply the Observation O2-1 and  obtain that for every  $N > 0$ 
$$
 \| V^N \|  \leq {c_0}^{2n_0+1} 
$$
 Well known
$$
 r(V) \leq  \| V^N \|^{1/N} \,, \qquad (N=1,2,\ldots) 
$$
 Hence
$$
 r(V)  \leq {c_0}^{(2n_0+1)/N} \,, \qquad (N=1,2,\ldots) 
$$
 and
$$
 r(V)  \leq 1 \,. 
$$
\\
$\Box$

\par\addvspace{2\medskipamount}\par\noindent
\par\noindent 
% %

\newpage 
\par\addvspace{2\medskipamount}\par\noindent
\par\noindent 
% %
{\bf Lemma 2.}
\par\noindent 
Let
\begin{eqnarray*}
|v_n|  
&\not=&
0
\end{eqnarray*}
and $V$ be a bilateral weighted shift defined by 
\begin{eqnarray*}
Vb_n  
&:=&
\frac{v_{n+1}}{v_n} b_{n+1}
\end{eqnarray*}
Suppose that
$$
     \exists c > 0  \exists n_0 \geq 0   
$$
so that
$$
  \forall n \qquad 
 |n| \geq n_0  \Rightarrow    \Bigl| \frac{v_{n+1}}{v_n} \Bigr| \leq  c   
\qquad\qquad\qquad\qquad(*)
$$
 Then
$$||V||<\infty$$
and 
\begin{eqnarray*}
 r(V)
&\leq&
c 
\end{eqnarray*}
As a consequence 
\begin{eqnarray*}
 r(V^*)=r(V)
&\leq&
c 
\end{eqnarray*}
\medskip\noindent
{\bf Proof}
\\
 Put 
$$
     w_{n} = \frac{v_n}{c^{n}} 
$$
 Then 
$$
   \Bigl| \frac{w_{n+1}}{w_n}\Bigr|
 = \Bigl| \frac{v_{n+1}}{v_n}\Bigr| \frac{1}{c}
$$
  Hence 
$$
   \Bigl| \frac{w_{n+1}}{w_n}\Bigr|
     \leq  1  \qquad (  |n| \geq |n_0| ) 
$$
Let $W$ be a bilateral weighted shift defined by 
\begin{eqnarray*}
Wb_n  
&:=&
\frac{w_{n+1}}{w_n} b_{n+1}
\end{eqnarray*}
By the Lemma 1 
$$
    \| W \| < \infty \,,\, \qquad r(W) \leq 1 \,.
$$
On the other hand 
$$
    W = \frac{1}{c} V \,.
$$
Hence
$$
    V = c W \,.
$$
Hence
$$
    r(V) = c \; r(W) \leq c \,.
$$
\\
$\Box$
 
\par\addvspace{2\medskipamount}\noindent 
% %

\par\noindent 
% %

%\input unitary-2013fev20-2.tex

\newpage

\newpage 
\par\noindent 
% %
{\bf Lemma 3}
\par\noindent 
Let
\begin{eqnarray*}
|v_n|  
&\not=&
0
\end{eqnarray*}
and $V$ be a bilateral weighted shift defined by 
\begin{eqnarray*}
Vb_n  
&:=&
\frac{v_{n+1}}{v_n} b_{n+1}
\end{eqnarray*}
Suppose that
$$
      \exists n_1 \geq 0  
$$
so that
$$
  \forall n  \qquad 
|n| \geq n_0  \Rightarrow    \Bigl| \frac{v_{n+1}}{v_n} \Bigr| \geq  1   
\qquad\qquad\qquad\qquad(*)
$$
 Then
$$||V^{-1}||<\infty$$
and 
\begin{eqnarray*}
 r(V^{-1})
&\leq&
1 
\end{eqnarray*}
As a consequence 
\begin{eqnarray*}
 r({V^{-1}}^*) = r(V^{-1})
&\leq&
1 
\end{eqnarray*}
\medskip\noindent
{\bf Proof}
\\
By supposition,
$$
|n| \geq n_1  \Rightarrow    \Bigl| \frac{v_{n+1}}{v_n} \Bigr| \geq  1   
%\qquad\qquad\qquad\qquad(*)
$$
 Hence 
$$
|n| \geq n_1  \Rightarrow    \Bigl| \frac{v_n}{v_{n+1}} \Bigr| \leq  1   
%\qquad\qquad\qquad\qquad(*)
$$
 and then 
$$
|n-1| \geq n_1  \Rightarrow    \Bigl| \frac{v_{n-1}}{v_{n}} \Bigr| \leq  1   
%\qquad\qquad\qquad\qquad(*)
$$
$$
|n| \geq n_1+1  \Rightarrow    \Bigl| \frac{v_{n-1}}{v_{n}} \Bigr| \leq  1   
%\qquad\qquad\qquad\qquad(*)
$$
\\
Put 
$$
  c_1 = max\{ 1, max \Bigl\{ \Bigl| \frac{v_{n-1}}{v_n} \Bigr| \quad 
               \Big| \quad n=-n_1-1, -n_1, \cdots , n_1+1  \Bigr\}\}
$$
Hence for every $n \in {\bf Z}$, $N > 0$ 
$$
  \Bigl| \frac{v_{n-N}}{v_n} \Bigr| = \Bigl| \frac{v_{n-1}}{v_n} \Bigr|
  \cdots
  \Bigl| \frac{v_{n-N}}{v_{n-(N-1)}} \Bigr|  \leq {c_1}^{2n_1+3} 
$$
Apply the Observation O2-1 and  obtain that for every  $N > 0$ 
$$
 \| V^{-N} \|  \leq {c_1}^{2n_1+3} 
$$
 Well known
$$
 r(V^{-1}) \leq  \| V^{-N} \|^{1/N} \,, \qquad (N=1,2,\ldots) 
$$
 Hence
$$
 r(V^{-1})  \leq {c_1}^{(2n_1+3)/N} \,, \qquad (N=1,2,\ldots) 
$$
 and
$$
 r(V^{-1})  \leq 1 \,. 
$$
\\
$\Box$

\par\addvspace{2\medskipamount}\par\noindent
\par\noindent 
% %

\par\addvspace{2\medskipamount}\par\noindent
\par\noindent

\newpage 
\par\noindent 
% %
{\bf Lemma 4}
\par\noindent 
Let
\begin{eqnarray*}
|v_n|  
&\not=&
0
\end{eqnarray*}
and $V$ be a bilateral weighted shift defined by 
\begin{eqnarray*}
Vb_n  
&:=&
\frac{v_{n+1}}{v_n} b_{n+1}
\end{eqnarray*}
Suppose
$$
     \exists c \geq 0   \exists n_1 \geq 0   
$$
so that
$$
  \forall n \qquad
 |n| \geq n_1  \Rightarrow    \Bigl|\frac{v_{n+1}}{v_n}\Bigr| \geq  c^{-1}   
\qquad\qquad\qquad\qquad(**)
$$
 Then
$$||V^{-1}||<\infty$$
and 
\begin{eqnarray*}
 r(V^{-1})
&\leq&
c 
\end{eqnarray*}
As a consequence 
\begin{eqnarray*}
 r({V^{-1}}^*) = r(V^{-1})
&\leq&
c 
\end{eqnarray*}
\medskip\noindent
{\bf Proof.}
\\
 Put 
$$
     w_{n} = {v_n}c^{n} 
$$
 Then 
$$
   \Bigl| \frac{w_{n+1}}{w_n}\Bigr|
 = \Bigl| \frac{v_{n+1}}{v_n}\Bigr| c 
$$
  Hence 
$$
   \Bigl| \frac{w_{n+1}}{w_n}\Bigr|
     \geq  1  \qquad (  |n| \geq |n_1| ) 
$$
Let $W$ be a bilateral weighted shift defined by 
\begin{eqnarray*}
Wb_n  
&:=&
\frac{w_{n+1}}{w_n} b_{n+1}
\end{eqnarray*}
By the Lemma 3 
$$
    \| W^{-1} \| < \infty \,,\, \qquad r(W^{-1}) \leq 1 \,.
$$
On the other hand 
$$
    W = c V \,.
$$
$$
    V^{-1} = c W^{-1} \,.
$$
Hence
$$
    r(V^{-1}) = c \; r(W^{-1}) \leq c \,.
$$
\\
$\Box$
 
\par\addvspace{2\medskipamount}\noindent

\newpage

\par\addvspace{2\medskipamount}\noindent 
{\bf Lemma 5}
\par\noindent 
Let
\begin{eqnarray*}
|v_n|  
&\not=&
0
\end{eqnarray*}
\begin{eqnarray*}
v_n  
&=&
v_{-n}
\end{eqnarray*}
and $V$ be a bilateral weighted shift defined by 
\begin{eqnarray*}
Vb_n  
&:=&
\frac{v_{n+1}}{v_n} b_{n+1}
\end{eqnarray*}
Suppose that
$$
     \exists c \geq1 \exists n_0 \geq 0   
$$
so that
$$
 \forall n \qquad  
|n| \geq n_0  \Rightarrow    \Bigl| \frac{v_{n+1}}{v_n} \Bigr| \leq  c    
\qquad\qquad\qquad\qquad(*)
$$
 Then
$$||V^{-1}||=||V||<\infty$$
and 
\begin{eqnarray*}
 r(V^{-1})=r(V)
&\leq&
c 
\end{eqnarray*}
As a consequence 
\begin{eqnarray*}
 r(V^*)=r(V)=r(V^{-1})=r({V^{-1}}^*)
&\leq&
c 
\end{eqnarray*}
\medskip\noindent
{\bf Proof}
\\
  Define $R$ by 
\begin{eqnarray*}
Rb_n  
&:=&
 b_{-n}
\end{eqnarray*}
   Then
$$
    R^{-1} = R , \qquad  R^2=I , \qquad ||R||=1<\infty , \qquad ||R^{-1}||=1<\infty  
$$
  Consider
$R^{-1}SR$
$$
     R^{-1}VRb_n=R^{-1}Vb_{-n}
    =  R^{-1}\frac{v_{-n+1}}{v_{-n}} b_{-n+1}
    =  R^{-1}\frac{v_{n-1}}{v_{n}} b_{-n+1}
    =\frac{v_{n-1}}{v_{n}} b_{n-1}
    = V^{-1}b_n
$$
   Hence
$$
       R^{-1}VR = V^{-1}
$$
\\
Now apply the Lemma 2.

\par\addvspace{2\medskipamount}\par\noindent

\newpage

\par\addvspace{2\medskipamount}\par\noindent
{\bf Lemma 6 }
\par\noindent 
Let
\begin{eqnarray*}
c  
& \geq &
1
\end{eqnarray*}
\begin{eqnarray*}
\phi(x)  
&:=&
x\sin(\frac{\pi}{2}\log_2(1+\log_2(1+x))) , \qquad ( x \geq 0)
\\
\psi(x)  
&:=&
x^{1/2}\sin(\frac{\pi}{2}\log_2(1+\log_2(1+x))) , \qquad ( x \geq 0)
\end{eqnarray*}
\begin{eqnarray*}
v_n  
&:=&
c^{ \phi(|n|)} e^{\psi{|n|}}
\end{eqnarray*}
so that 
\begin{eqnarray*}
v_n  
& = &
\overline{ v_n } 
\end{eqnarray*}
\begin{eqnarray*}
v_n  
& > &
0 
\end{eqnarray*}
\begin{eqnarray*}
v_n  
&=&
v_{-n}
\end{eqnarray*}
and $V$ be a bilateral weighted shift defined by 
\begin{eqnarray*}
Vb_n  
&:=&
\frac{v_{n+1}}{v_n} b_{n+1}
\end{eqnarray*}
Then 
\begin{eqnarray*}
 r(V^*)=r(V)=r(V^{-1})=r({V^{-1}}^*)
& = &
c 
\end{eqnarray*}
and 
$$
S(V,c)=\{0\} \,,\,
S(V^{-1},c)=\{0\} \,,\,
S(V^{*-1},c)=\{0\} \,,\,
S(V^*,c)=\{0\} \,.\,
$$
\medskip\noindent
{\bf Proof}
\begin{eqnarray*}
z
&:=&
\frac{\pi}{2}\log_2(1+\log_2(1+x))
\end{eqnarray*}
and note 
\begin{eqnarray*}
\phi'(x)
&=&
\sin(z) 
+ x\cos(z )\frac{\pi}{2}
\frac{1}{1+\log_2(1+x)}
\frac{1}{1+x}\frac{1}{(ln2)^2}
 \\
&=&
\sin(z) 
+ \cos(z )
\frac{1}{1+\log_2(1+x)}
\frac{x}{1+x} k , \qquad k=\frac{\pi}{2(ln2)^2}
\end{eqnarray*}
   Hence
\begin{eqnarray*}
\phi'(0)
&=&
0
\end{eqnarray*}
   and if  
\begin{eqnarray*}
 x
&\geq&
  2^{2k/\epsilon} -1 
\end{eqnarray*}
   then 
\begin{eqnarray*}
|\phi'(x)|
&\leq&
 1 +\frac{1}{(1+2k/ \epsilon)}k
\\
&\leq&
 1 +\epsilon/2
\\
\end{eqnarray*}
   Hence if  
\begin{eqnarray*}
 n
&\geq&
  2^{2k/ \epsilon} -1 
\end{eqnarray*}
   then 
\begin{eqnarray*}
|\phi(|n+1|)-\phi(|n|)|
&\leq&
 1 +\epsilon/2
\\
\end{eqnarray*}
 Now note that if $ x\not=0$ then 
\begin{eqnarray*}
\psi'(x)
&=&
\frac{1}{2x^{1/2}}\sin(z) 
+ x^{1/2}\cos(z )\frac{\pi}{2}
\frac{1}{1+\log_2(1+x)}
\frac{1}{1+x}\frac{1}{(ln2)^2}
 \\
&=&
\frac{1}{2x^{1/2}}\sin(z) 
+ \cos(z )
\frac{1}{1+\log_2(1+x)}
\frac{x^{1/2}}{1+x} k , \qquad k=\frac{\pi}{2(ln2)^2} . 
\end{eqnarray*}
 As a result, 
$$
 \psi'(x) \to 0 \mbox{ as } x \to \pm \infty
$$
 and  
$$
 \psi(|n+1|)-\psi(|n|) \to 0 \mbox{ as } n \to \pm \infty
$$
  Hence if 
$$ 
   c> 1
$$
   then 
$$
   \forall  \epsilon>0  \exists n_0 \geq 0   
$$
so that
$$
\forall n \qquad 
   |n| \geq n_0  \Rightarrow    \Bigl| \frac{v_{n+1}}{v_n} \Bigr| 
   \leq  c^{1+\epsilon} 
 $$
 Then apply the Lemma 5 and obtain
$$||V^{-1}||=||V||<\infty$$
and 
$$ \forall \epsilon > 0 $$
\begin{eqnarray*}
  r(V^*)=r(V)=r(V^{-1})=r({V^{-1}}^*)
&\leq&
c^{1+\epsilon} 
\end{eqnarray*}
Hence 
\begin{eqnarray*}
  r(V^*)=r(V)=r(V^{-1})=r({V^{-1}}^*)
&\leq&
c 
\end{eqnarray*}
 The case $c=1$ is obvious.

\par\addvspace{2\medskipamount}\par\noindent
 Now choose two sequences of integers defining them by 
$$
 n_k\,:=\,2^{2^{1+4k}-1}-1;\quad m_k\,:=2^{2^{3+4k}-1}-1  \qquad (k=1,2,...) \,.
$$
 Then 
$n_k,\, m_k \in {\bf N}$,\ 
$n_k\,\to\,+\infty$,$m_k\,\to\,+\infty$
 (as 
$k\,\to\,+\infty$),
 and simultaneously 
$$
 v_{n_k} = v_{-n_k} = c^{n_k}e^{\sqrt{n_k}} \,; \quad
 v_{m_k}^{-1} = v_{-m_k}^{-1} = c^{m_k}e^{\sqrt{m_k}} \,.
$$

\noindent
 We see that no estimation of the form 

$$
\begin{array}{clc}
 |v_N|\,\leq\,M'c^N,\,
&
 |v_{-N}|\,\leq\,M'c^N,\,
&
 |v_{N}|^{-1}\,\leq\,M'c^N,\,
\\
&
 |v_{-N}|^{-1}\,\leq\,M'c^N
&
\end{array}
$$

\noindent 
 (for 
$N=0,1,\cdots$
 )
 is possible. 
 On looking at the Observation {\bf O2-2}, we see that 
$$
S(V,c)=\{0\} \,,\,
S(V^{-1},c)=\{0\} \,,\,
S(V^{*-1},c)=\{0\} \,,\,
S(V^*,c)=\{0\} \,.\,
$$
 This is just what was to be proven.
\\
$\Box$

\par\addvspace{2\medskipamount}\par\noindent
\par\noindent 
% %

%\input unitary-2013fev20-5.tex

%
%
%
%
%
%
%
\newpage

%{\bf 2013 mart 09}
%
\par\addvspace{2\medskipamount}\noindent 
{\bf Lemma 7 }
\begin{eqnarray*}
\{ \widehat{U}^N(b_0\oplus b_0), \widehat{U}^M(b_0\oplus b_0) \}
&=& 
0 , \quad (M \not = N)
\\
\{ \widehat{U}^N(b_0\oplus b_0), \widehat{U}^N(b_0\oplus b_0) \}
&=& 
2(b_0,b_0) = 2 > 0
\\
\{ \widehat{U}^N(b_0\oplus -b_0), \widehat{U}^M(b_0\oplus -b_0) \}
&=&
0 , \quad (M \not = N)
\\
\{ \widehat{U}^N(b_0\oplus -b_0), \widehat{U}^N(b_0\oplus -b_0) \}
&=& 
-2(b_0,b_0) = -2 < 0
\\
\{ \widehat{U}^N(b_0\oplus b_0), \widehat{U}^M(b_0\oplus -b_0) \}
&=&
0 .
\\
\end{eqnarray*}

\par\addvspace{2\medskipamount}\noindent 
{\bf Definition}
\begin{eqnarray*}
 L_{+}
& := &
 span\{\widehat{U}^N(b_0\oplus b_0)   | N \in {\bf Z}\}
 \equiv span\{(U^N b_0\oplus {U^*}^{-N} b_0)   | N \in {\bf Z}\}
\\
 L_{-}
& := &
 span\{\widehat{U}^N(b_0\oplus -b_0)   | N \in {\bf Z}\}
 \equiv span\{(U^N b_0\oplus -{U^*}^{-N} b_0)   | N \in {\bf Z}\}
\\
\end{eqnarray*}
 Then  
\begin{eqnarray*}
(1)
&&
  \widehat{U}^{\pm 1}L_{+} = L_{+}
\\
(2)
&&
  \widehat{U}^{\pm 1}L_{-} = L_{-}
%\\
%(3)
%&&
%   L_{+} + L_{-} \supset \{ b_N \oplus 0 | N \in {\bf Z}\}
%                    \cup \{ 0 \oplus b_N | N \in {\bf Z}\}   
\\
(3)
&&
b_N \oplus 0 \in L_{+} + L_{-}   \quad ( N \in {\bf Z} )
\\
%(3)
&&
0 \oplus b_N \in L_{+} + L_{-}  \quad ( N \in {\bf Z} )
\\
(4)
&&
  \mbox{if } \{ L_{+} + L_{-}, x \}=0 \mbox{ then } x = 0 
\\
(5)
&&
  \overline{ L_{+} + L_{-} }= \widehat{ H  }
\\
(6)
&&
  \{ L_{+} , L_{-}\}=\{0\} 
\\
(7)
&&
 \{x,x\} > 0 \quad ( x \in L_{+}\backslash \{0\} ) 
\\
(8)
&&
 \{x,x\} < 0 \quad ( x \in L_{-}\backslash \{0\} ) 
\\
(9)
&&
  L_{+} \cap  L_{-} = \{0\}
\\
\end{eqnarray*}
 Moreover,   
\begin{eqnarray*}
(6') 
&&
 \{ \overline{L_{+}} , \overline{L_{-}}\}=\{0\} 
\\
(7')
&&
 \{x,x\} > 0 \quad ( x \in \overline{L_{+}}\backslash \{0\} ) 
\\
(8')
&&
 \{x,x\} < 0 \quad ( x \in \overline{L_{-}}\backslash \{0\} ) 
\\
(9')
&&
 \overline{L_{+}} \cap  \overline{L_{-}} = \{0\}
\\
\end{eqnarray*}
 If in addition $U$ and $U^{-1}$ are bounded, then 
\begin{eqnarray*}
(1')
&&
  \widehat{U}^{\pm 1}\overline{L_{+}} = \overline{L_{+}}
\\
(2')
&&
  \widehat{U}^{\pm 1}\overline{L_{-}} = \overline{L_{-}}
\end{eqnarray*}
\par\addvspace{2\medskipamount}\noindent 
{\bf Proof} is straightforward.

\par\addvspace{2\medskipamount}\noindent 
 {\bf Remark}
\par\noindent 
 $$ g\oplus h \in L_{+} $$
 $$ g\oplus h
     = \sum_n g(n)b_n \oplus \sum_n \frac{|u_0|^2}{|u_n|^2}g(n)b_n $$
 $$ g\oplus h
     = g \oplus \sum_n \frac{|u_0|^2}{|u_n|^2}b_n(b_n,g) $$
 $$ g\oplus h \in L_{-} $$
 $$ g\oplus h
     = \sum_n g(n)b_n \oplus -\sum_n \frac{|u_0|^2}{|u_n|^2}g(n)b_n $$
 $$ g\oplus h
     = g \oplus -\sum_n \frac{|u_0|^2}{|u_n|^2}b_n(b_n,g) $$
$\Box$ 

 $$ g\oplus h \in L_{+} $$
 $$ g\oplus h
     = \sum_n g(n)b_n \oplus \sum_n \frac{|u_0|^2}{|u_n|^2}g(n)b_n $$

\begin{eqnarray*}
(h, U^N g)
&=&
  (\sum_n \frac{|u_0|^2}{|u_n|^2}g(n)b_n , U^N \sum_m g(m)b_m ) 
\\
&=&
(\sum_n \frac{|u_0|^2}{|u_n|^2}g(n)b_n ,
\sum_m \frac{u_{m+N}}{u_m}g(m)b_{m+N} )
\\
&=&
\sum_m \frac{|u_0|^2}{|u_{m+N}|^2}\overline{g(m+N)} \frac{u_{m+N}}{u_m}g(m))
\\
&=&
\sum_m \frac{|u_0|^2}{\overline{u_{m+N}}}\overline{g(m+N)} \frac{1}{u_m}g(m))
\end{eqnarray*}

\par\addvspace{2\medskipamount}\noindent 
 {\bf Remark} to the Lemma 1 + Lemma 2, Lemma 3 + Lemma 4, Lemma 5. 
\par\noindent 
The statements look very naturally, and I think may be known. 
But I have not found.

\end{document}